\documentclass[twoside,11pt]{article}
\usepackage{graphicx,amscd,amsfonts,amsmath,amsthm,amssymb,latexsym,amsfonts,color}
\usepackage{epsfig,float,epstopdf,array,bm,cite}
\usepackage{cases}
\usepackage{array}
\usepackage{multirow}
\usepackage[table]{xcolor}
\usepackage{relsize}
\usepackage{graphicx}
\usepackage{float}
\usepackage[bookmarksnumbered, colorlinks,plainpages]{hyperref}
\def\diag{\textrm{diag}}

\def\h{\hspace{-0.2cm}}

\newtheorem{theorem}{\bf Theorem}
\newtheorem{lemma}{\bf Lemma}

\newtheorem{example}{\bf Example}

\newtheorem{remark}{\bf Remark}
\begin{document}
\title{\bf On the preconditioning of three-by-three block saddle point problems}
\author{\small\bf  Hamed Aslani$^\dag$, Davod Khojasteh Salkuyeh$^{\dag,\ddag}$\thanks{\noindent Corresponding author. \newline
		Emails:	hamedaslani525@gmail.com (H. Aslani), khojasteh@guilan.ac.ir (D.K. Salkuyeh), f.beik@vru.ac.ir (F.P.A. Beik) },
	Fatemeh Panjeh Ali Beik$^\S$ \\[2mm]
	\textit{{\small $^\dag$Faculty of Mathematical Sciences, University of Guilan, Rasht, Iran}} \\
	\textit{{\small $^\ddag$Center of Excellence for Mathematical Modelling, Optimization and Combinational}}\\
	\textit{{\small Computing (MMOCC), University of Guilan, Rasht, Iran}} \\
	\textit{{\small $^\S$Department of Mathematics, Vali-e-Asr University of Rafsanjan, P.O. Box 518, Rafsanjan, Iran}} \\
}
\date{}
\maketitle
\vspace{-0.5cm}

\noindent\hrulefill\\
{\bf Abstract.}   We establish a new iterative method for solving a class of large and sparse linear systems of equations with three-by-three block coefficient matrices having saddle point structure. Convergence properties of the proposed method are studied in details and its induced preconditioner is examined for accelerating the convergence speed of generalized minimal residual (GMRES) method. More precisely, we analyze the eigenvalue distribution of the preconditioned  matrix. Numerical experiments are reported to demonstrate the effectiveness of the proposed preconditioner. \\[-3mm]


\noindent\hrulefill\\

\pagestyle{myheadings}\markboth{H. Aslani, D.K. Salkuyeh, F.P.A. Beik  }{On block iterative scheme and preconditioner for the  3$\times$3 block saddle point problem}
\thispagestyle{empty}

\section{Introduction}
Consider the following three-by-three block system of linear equations,
\begin{equation}\label{eq1}
\mathcal{A} {\bf x} \equiv\left(\begin{array}{ccc}
{A} & {B^{T}} & {0} \\
{B} & {0} & {C^{T}} \\
{0} & {C} & {0}
\end{array}\right)\left(\begin{array}{l}
{x} \\
{y} \\
{z}
\end{array}\right)=\left(\begin{array}{l}
{f} \\
{g} \\
{h}
\end{array}\right),
\end{equation}
where $A\in \mathbb{R}^{n\times n}$, $B\in \mathbb{R}^{m\times n}$, $C\in \mathbb{R}^{l\times m}$, $f\in \mathbb{R}^n$, $g\in \mathbb{R}^m$ and $h\in \mathbb{R}^l$ are known, and ${\bf x}=\left(x; y; z\right)$ is an unknown vector to be determined. Here, the \textsc{Matlab} symbol $(x;y;z)$ is utilized to denote the vector $(x^{T},y^{T},z^{T})^{T}.$

In the sequel, we assume that the matrix $A$ is a symmetric positive definite and the matrices $B$ and $C$ have full row rank. These assumptions guarantee the existence of a unique solution of \eqref{eq1}; see \cite{A2} for further details.

Evidently matrix $\cal A$ can be regarded as a $2\times 2$ block matrix using
the following partitioning strategy,
\begin{equation}\label{part}
\mathcal{A} = \left( {\begin{array}{cc|c}
	A & {B^T } & {0}  \\
	B & 0 & C^T   \\
		\hline
	0 & C & 0  \\
	\end{array}} \right).
\end{equation}
As seen, the above block matrix has a saddle point structure. Hence, we call Eq. \eqref{eq1} by three-by-three block saddle point problem.

Linear system of the form \eqref{eq1}  arises from many practical scientific and engineering application backgrounds, e.g., the discrete finite element methods for solving time-dependent Maxwell equation with discontinuous coefficient  \cite{Assous,A4,A5,A6},  the least squares problems \cite{A7}, the  Karush-Kuhn-Tucker (KKT)  conditions of a type of quadratic program \cite{A8} and so on.
Since the matrices $A,$ $B$ and $C$ in \eqref{eq1} are large and sparse, the solution of \eqref{eq1} is suited by iterative methods. In practice, stationary iterative methods may converge too slowly or fail to converge. For this reason they are usually combined with acceleration schemes, like Krylov subspace methods \cite{A9}. Here, we focus on preconditioned Krylov subspace methods, especially, the preconditioned GMRES method.

As seen, the coefficient matrix $\mathcal{A}$ in Eq. \eqref{eq1} can be considered in a two-by-two block form given by \eqref{part}. The observation was used in the literature for constructing preconditioners to improve the convergence speed of Krylov subspace methods for solving \eqref{eq1}, such as block triangular preconditioners \cite{A11,A12,A13,Benzi1,Benzi2,Beik1,Beik2}, shift-splitting preconditioners \cite{A14} and parameterized preconditioners \cite{A15}; for more details see also \cite{Benzi1,BenziSIAM,Salkuyeh}. Recently, Huang and Ma \cite{A1} proposed the following block diagonal preconditioner,
\begin{equation}\label{eq992}
\mathcal{P}_{D}=\left(\begin{array}{ccc}
{A} & {0} & {0} \\
{0} & {S} & {0} \\
{0} & {0} & {C S^{-1} C^{T}}
\end{array}\right),
\end{equation}
for solving \eqref{eq1} in which $S=B A^{-1} B^{T}.$ They also derive all the eigenpairs of preconditioned matrix. Xie and Li \cite{A2} presented the following three preconditioners
\[
{\cal P}_1=\begin{pmatrix}
A &   0  & 0 \\
B &  -S  & C^T \\
0 &   0  & CS^{-1}C^T
\end{pmatrix},~
{\cal P}_2=\begin{pmatrix}
A &   0  & 0 \\
B &  -S  & C^T \\
0 &   0  & -CS^{-1}C^T
\end{pmatrix},~
{\cal P}_3=\begin{pmatrix}
A &   B^T  & 0 \\
B &  -S    & 0 \\
0 &   0    & -CS^{-1}C^T
\end{pmatrix},
\]
and analyzed spectral properties of corresponding preconditioned matrices in the case $S=BA^{-1}B^T$. The reported numerical results in \cite{A2} show that the above preconditioners can significantly improve the convergence speed of GMRES method. It can be observed that the preconditioner ${\cal P}_1$ outperforms other preconditioners in terms of both required CPU time and number of iterations for the convergence.

Here, we consider the following equivalent form of \eqref{eq1}:
\begin{equation}\label{eq111}
{\cal B}{\bf x}\equiv \left(\begin{array}{ccc}
{A} & {B^{T}} & {0} \\
-{B} & {0} & -{C^{T}} \\
{0} & {C} & {0}
\end{array}\right)\left(\begin{array}{l}
{x} \\
{y} \\
{z}
\end{array}\right)=\left(\begin{array}{l}
{f} \\
{-g} \\
{h}
\end{array}\right)=\bf{b}.
\end{equation}
 Although the coefficient matrix of the system \eqref{eq111} is not symmetric, it has some desirable properties. For instance, the matrix {$\mathcal{B}$} is positive semidefinite, i.e., $\mathcal{B}+\mathcal{B}^{T}$ is symmetric positive semidefinite. This is a significant for the GMRES method. In fact, the restarted version of GMRES($m$) converges for all $m\geq 1$. Recently, some iterative schemes have been extended in the literature for solving \eqref{eq111}. For instance, Cao \cite{CAOAML} presented the shift-splitting method. In \cite{Huang-NumerAlgor,Huang-NLWA}, the Uzawa-type methods were developed.  In this work, we present  a new type of iterative method for solving three-by-three block saddle point problem \eqref{eq111}. Next, we extract a preconditioner from the presented iterative method and examine its performance for speeding up the convergence of GMRES.

The remainder of this paper organized as follows. Before ending this section, we present notations and basic preliminaries used in next sections. In section \ref{sec2}, we propose a new iterative method for solving \eqref{eq111} and study its converges properties. In section \ref{sec3}, we extract a preconditioner from the proposed method and analyze the spectrum of preconditioned matrix. Brief discussions are given in section \ref{sec4} about practical implementation of the preconditioner. In section \ref{sec5}, we report some numerical results and  brief concluding remarks are included in section \ref{sec6}.

 {Throughout this paper, the identity matrix is denoted by  $I$. The symbol $x^{*}$ is used for the conjugate transpose of the vector $x.$
 	%For a given matrix $A\in \mathbb{R}^{n\times n}$, notations $A^{-1}$ and $A^{T}$ are used for the transpose and inverse of the matrix $A,$ respectively.
 	For any square matrix $A$ with real eigenvalues, the minimum and maximum eigenvalues of $A$ are indicated by $\lambda_{\min} (A)$ and $\lambda_{\max} (A)$, respectively. The notation $\rho(A)$ stands for the spectral radius of $A.$ The matrix  $A\in \mathbb{R}^{n \times n}$ is called symmetric positive definite (SPD), if $A^{T}=A$ and $x^{T} A x>0$ for all nonzero $x\in \mathbb{R}^{n}.$ Similarly,  the matrix $A$ is called symmetric positive semidefinite (SPSD), if $A^{T}=A$ and $x^{T}A x \geqslant 0$ for all $x\in \mathbb{R}^{n}.$
 	We write $A \succ 0$ $ (A \succeq 0 ),$ if $A$ is SPD (SPSD). For two given matrices $A$ and $B,$ $A\succ B$ $(A \succeq B)$ means that $A-B \succ 0$ $(A-B \succeq 0).$  The matrix $A\in \Bbb{R}^{n\times n}$ is said to be positive (semi-) definite, if $A+A^T$ symmetric positive (semi-) definite. For any matrix $W$, we shall write its null space as $null(W).$  The norm $\parallel . \parallel$ indicates the 2-norm. }

\section{The proposed iteration scheme}\label{sec2}
Let us first consider the following splitting for the coefficient matrix $\mathcal{B}$ in \eqref{eq111}:
\begin{equation}\label{Split}
\mathcal{B}=\mathcal{P}-\mathcal{R},
\end{equation}
where
\begin{equation*}
\mathcal{P}=\left(\begin{array}{ccc}
{A} & {B^{T}} & {0} \\
{0} & {S} & {-C^{T}} \\
{0} & {C} & {0}
\end{array}\right), \quad
\mathcal{R}=\left(\begin{array}{ccc}
{0} & {0} & {0} \\
{B} & {S} & {0} \\
{0} & {0} & {0}
\end{array}\right),
\end{equation*}
in which $S$ is a given symmetric positive definite matrix. It is not difficult to verify that the matrix ${\cal P}$ is nonsingular.
The iteration scheme associated with splitting \eqref{Split} is given by
\begin{equation}\label{qe3}
{\bf x}^{(k+1)}=\mathcal{G} {\bf x}^{(k)}+{\bf c}, \quad k=0,1,2, \ldots,
\end{equation}
where ${\bf x}^{(0)}$ is an initial guess, $\mathcal{G}=\mathcal{P}^{-1} \mathcal{R} $ is the iteration matrix and ${\bf c}=\mathcal{P}^{-1} {\bf b}$.

Now, we present sufficient conditions under which the iterative scheme \eqref{qe3} is convergent. To this end, we first need to recall the following theorem.
\begin{theorem}\label{thorn} \cite[Theorem 7.7.3]{Horn1}
	Let $A$ and $B$ be two $n\times n$ real symmetric matrices such that $A$
	is positive definite and $B$ is positive semidefinite. Then $A \succeq  B$
	if and only if $\rho(A^{-1}B) \le 1,$ and $A\succ B$ if and only if $\rho(A^{-1}B) < 1$.\\
\end{theorem}

\begin{theorem}\label{th1}
	Let $A\succ 0$, $S\succ 0$ and $B$ and $C$ be  full row rank matrices. If
	$2S \succ B A^{-1} B^{T}$
	then the iterative method \eqref{qe3} converges to the unique solution of \eqref{eq111} for any initial guess.
	\begin{proof}
		Let $\lambda $ be an arbitrary eigenvalue of $\mathcal{G}=\mathcal{P}^{-1} \mathcal{R}$ with the corresponding eigenvector $w=\left(x; y; z\right)$. Consequently, we have $\mathcal{R} w=\lambda \mathcal{P}w$ which is equivalent to say that
		\begin{numcases}{}
		\lambda (A x+B^{T} y) =0,  \label{eq55} \\
	    \lambda ( S y - C^{T} z) = 	B x+ S y, \label{eq56}
		\\\lambda C y =0 .\label{eq57}
		\end{numcases}
		Without loss of generality, we may assume that $\lambda \neq 0$. Obviously $y \neq 0$, otherwise in view of the positive definiteness of $A$ and the assumption that $C$ has full row rank we conclude that $x$ and $z$ are both zero vectors which is in contradiction with the
		fact that  $(x;y;z)$ is an eigenvector. From Eqs. \eqref{eq55} and \eqref{eq57} we can deduce that
		$$x=-A^{-1} B^{T} y,\quad y^{*} C^{T}= 0.$$
		Multiplying both sides of Eq.  \eqref{eq56} on the left by $y^{*}$ and  substituting the preceding equalities, we derive
		\[
		\lambda =1 - \frac{y^{*} B A^{-1} B^{T} y}{y^{*} S y}.
		\]
		This shows that the eigenvalues of ${\cal G}$ are all real. By Theorem \ref{thorn}, it is immediate to conclude that $\lambda_{\max}(S^{-1}BA^{-1}B^T)=\rho(S^{-1}BA^{-1}B^T)<2$ if and only if $2S \succ B A^{-1} B^{T}$. This fact together with Courant-Fisher inequality \cite{A9} can deduce that
		\[
		0 < \frac{y^{*} B A^{-1} B^{T} y}{y^{*} S y} \le \lambda_{\max}(S^{-1}BA^{-1}B^T) <2.
		\]
		Therefore, we have
		\[
         |1 - \frac{y^{*} B A^{-1} B^{T} y}{y^{*} S y}|<1,
		\]
		which completes the proof.
	\end{proof}
\end{theorem}

We complete this section with a remark providing alternative sufficient conditions for convergence of iterative method  \eqref{qe3} which are stronger than  $2S \succ B A^{-1} B^{T}$, however, it might be easier to check the following sufficient conditions in some cases. To do so, we first remind the following two lemmas. The first one is a consequence of Weyl's
Theorem, see \cite[Theorem 4.3.1]{Horn1}.

\begin{lemma}\label{lem2.8} %\cite{Horn}
	Suppose that $A$ and $B$ are two Hermitian matrices. Then,
	\begin{eqnarray*}
		\lambda_{\max}(A+B) & \leq & \lambda_{\max}(A) + \lambda_{\max}(B),\\
		\lambda_{\min}(A+B) & \geq & \lambda_{\min}(A) + \lambda_{\min}(B).
	\end{eqnarray*}
\end{lemma}

\begin{lemma}\label{lem2} \cite{Zhang}
	Suppose that $A$ is a Hermitian negative definite matrix and
	$B$ is Hermitian positive semidefinite. Then the eigenvalues of $AB$
	are real and satisfy
	\[{\lambda _{\min}}(A){\lambda _{\min}}(B) \le {\lambda _{\max}}(AB)
	\le {\lambda _{\max}}(A){\lambda _{\min}}(B),\]
	\[{\lambda _{\min}}(A){\lambda _{\max}}(B) \le {\lambda _{\min}}(AB)
	\le {\lambda _{\max}}(A){\lambda _{\max}}(B).\]
\end{lemma}

\begin{remark}
Notice that  $2S \succ B A^{-1} B^{T}$ is equivalent to say that all eigenvalues of
$2S - B A^{-1} B^{T}$ are positive, i.e., $\lambda_{\min}(2S - B A^{-1} B^{T}) >0$.
From Lemma \ref{lem2.8}, it can be seen that
	\begin{align}
	\lambda_{\max }\left (B A^{-1} B^{T}\right)<2 \lambda_{\min}(S) , \label{qe4}
	\end{align}
	implies that $2S \succ B A^{-1} B^{T}$.
	Using Lemma \ref{lem2}, one can deduce that the condition \eqref{qe4} is satisfied as soon as
	\[
	\|B\|^2 < 2\lambda_{\min }(A)\lambda_{\min }(S),
	\]
	which follows from the fact that $$	\lambda_{\max }\left (B A^{-1} B^{T}\right)=	\lambda_{\max }\left ( A^{-1} B^{T}B\right) \le \lambda_{\max }\left ( A^{-1} \right)\lambda_{\max }\left (B^{T}B\right)=\frac{\|B\|^2}{\lambda_{\min}(A)}.$$
\end{remark}

%\begin{remark}
%	Assume that the conditions of Theorem \ref{th1} are satisfied.
%	If
%	$S\succeq B A^{-1} B^{T}$, then the iterative method \eqref{qe3} for solving the saddle point problem \eqref{eq1} is unconditionally convergent. In the special case, if $S=B A^{-1} B^{T}$ then $\rho({\cal G})=0$.
%\end{remark}

\section{The induced preconditioner and its spectral analysis}\label{sec3}

     From the splitting \eqref{Split} we have
     \[
     {\cal P}^{-1}{\cal B}=I-{\cal P}^{-1}{\cal R}=I-{\cal G}.
     \]
     Therefore, under the conditions of Theorem \ref{th1} the eigenvalues of ${\cal P}^{-1}{\cal B}$ are contained in the interval  $(0,2]$.
     Thus,
     \begin{equation}\label{eq99}
     \mathcal{P}=\left(\begin{array}{ccc}
     {A} & {B^{T}} & {0} \\
     {0} & {S} & {-C^{T}} \\
     {0} & {C} & {0}
     \end{array}\right),
     \end{equation}
     can be used as a preconditioner to accelerate the convergence of Krylov subspace methods like GMRES for solving the system \eqref{eq111}.

     In the succeeding theorem, we investigate the spectral properties of ${\cal P}^{-1}{\cal B}$ in more details.

     %\subsection{Spectral analysis of ${\cal P}^{-1}{\cal B}$}

     \begin{theorem}\label{th5}
     	Let $A$ be symmetric positive definite and $B$ and $C$ be of full row rank.
     	Then all the eigenvalues of $\mathcal{P}^{-1} \mathcal{B}$ are real and nonzero. Furthermore, $\lambda =1 $ is an eigenvalue of algebraic multiplicity at least $n+l$ and its corresponding eigenvectors are of the form
     	 $(x;-S^{-1} B x;z)$ where $x\in\Bbb{R}^{n}$ and $ z \in \Bbb{C}^l$ such that $x,z$ are not simultaneously zero.\\
     	The remaining eigenvalues of $\mathcal{P}^{-1} \mathcal{B}$  are of the form
     	\[
     	\lambda=\frac{y^{*} B A^{-1} B^{T} y}{y^{*} S y},
     	\]
     	and the corresponding eigenvectors are of the form $(-A^{-1} B^{T}y;y;z)$ for all $0\neq y \in null(C)$ and arbitrary $z$.
     	\begin{proof}
     		Let $\lambda$ be an arbitrary eigenvalue of  $\mathcal{P} ^{-1} \mathcal{B}$ with the corresponding eigenvector $(x;y;z)$, i.e.,
     		\begin{numcases}{}
     		A x+B^{T} y=\lambda (A x + B^{T} y), \label{eq22}\\ -B x - C^{T} z=\lambda ( S y - C^{T} z),
     		\label{eq23}
     		\\C y =\lambda Cy .\label{eq24}
     		\end{numcases}
     		Let $x=0$. If $\lambda \neq 1,$ then by \eqref{eq22} we have $B^{T}y=0,$ which shows taht $y=0.$ This
     		along with \eqref{eq23} leads to $C^{T}z = 0.$ Since $C$ is a full row rank matrix, then $z=0$. Consequently, we have $(x; y; z) = (0; 0; 0)$ which contradicts
     		with the fact that $(x;y;z)$ is an eigenvector. If $\lambda = 1$, then by \eqref{eq23} and the positive definiteness of $S$ we derive that $y=0.$ In addition, the corresponding eigenvectors are $(0;0;z),$ with $z \neq 0.$ In fact, $\lambda=1$ is an eigenvalue of $\mathcal{P}^{-1} \mathcal{B}$ with multiplicity $l$ corresponding eigenvector $(0;0;z)$ with $0\neq z \in \Bbb{C}^l$.
     		
     		In the following, we consider the case that $x \neq 0.$ If $y=0$, then Eqs.
     		\eqref{eq22} and \eqref{eq23} are reduced to
     		\begin{align}
     		A x=\lambda A x \quad \text{and} \quad -Bx -C^{T} z=-\lambda C^{T}z, \label{eq25}
     		\end{align}
     		respectively. The first relation shows that $ \lambda = 1.$ By substituting it into the second equality of \eqref{eq25}, we have $Bx=0.$ Therefore, the corresponding eigenvectors are of the form $(x;0;z)$ with $0 \ne x \in null(B)$ and $z\in\Bbb{R}^{l}$. Notice that, in general, we can observe that $\lambda=1$ and $(x;0;z)$ is an eigenpair of $\mathcal{P}^{-1} \mathcal{B}$ where $x \in null(B)$ and $x,z$ are not simultaneously zero.
     		
     		In summary,  using \eqref{eq23} and in view of the positive definiteness of $S$, we can conclude that  $\lambda = 1 $ and $(x;-S^{-1} B x;z)$
     		is an eigenpair of $\mathcal{P}^{-1} \mathcal{B}$.
     		
     		It is immediate to see that if $x$ and $y$ are both zero vectors then  $\lambda=1$ and $z$ must be a nonzero vector. In rest of the proof, we assume that $x \neq 0$ and $y \neq 0.$ If $\lambda \neq 1$, then from  Eqs. \eqref{eq22} and \eqref{eq24}, we observe that $x=-A^{-1} B^{T}y$ and $Cy=0$, respectively. Pre-multiplying both sides of \eqref{eq23} from left by $y^{*}$ and substituting deduced $x$ and $z$ into \eqref{eq23}, we get
     		$$\lambda=-\frac{y^{*} Bx}{y^{*} S y}=\frac{y^{*} B A^{-1} B^{T} y}{y^{*} S y}.$$
     		Hence, the corresponding eigenvectors are of the form $(-A^{-1} B^{T}y;y;z)$ for all $0\neq y \in null(C)\subseteq \Bbb{R}^{m}$ and arbitrary $z$.
     	\end{proof}
     \end{theorem}

\begin{remark}
Let $S$ be an arbitrary symmetric positive definite matrix. From Theorem \ref{th5} we see that the non-unit eigenvalues of the preconditioned matrix $\mathcal{P}^{-1} \mathcal{B}$ satisfies
\[
0<\frac{\lambda_{\min}(B A^{-1} B^{T})}{\lambda_{\max}(S)}\leq \lambda=\frac{y^{*} B A^{-1} B^{T} y}{y^{*} S y} \leq \frac{\lambda_{\max}(B A^{-1} B^{T})}{\lambda_{\min}(S)}.   	
\]
\end{remark}

\begin{theorem}\label{minpol}
Under the assumptions of Theorem \ref{th5}, if $S=B A^{-1} B^{T},$ then all  the eigenvalues of preconditioned matrix
${\cal H}={\cal P}^{-1}{\cal B}$ are equal to $1$ and its minimal polynomial is of degree 2.
\begin{proof}
Consider the matrix $\mathcal{P}_{D}$ defined in Eq. \eqref{eq992} with	$S=B A^{-1} B^{T}$. Obviously, $\mathcal{P}_{D}$ is symmetric positive definite, therefore there is a symmetric positive definite matrix   $\mathcal{P}_{D}^\frac{1}{2}$ such that $\mathcal{P}_{D}=\mathcal{P}_{D}^\frac{1}{2}\mathcal{P}_{D}^\frac{1}{2}$. Similar to the proof of Theorem 3.1 in \cite{A2}, we see that the matrix ${\cal H}$ is similar to the matrix
\begin{eqnarray}
\hat{{\cal H}}&:=& \mathcal{P}_{D}^\frac{1}{2} {\cal H} \mathcal{P}_{D}^{-\frac{1}{2}}\\
&=&    \mathcal{P}_{D}^\frac{1}{2} {\cal P}^{-1}{\cal B} \mathcal{P}_{D}^{-\frac{1}{2}}  \\
&=&  \left(\mathcal{P}_{D}^{-\frac{1}{2}}{\cal P} \mathcal{P}_{D}^{-\frac{1}{2}}\right)^{-1}  \left( \mathcal{P}_{D}^{-\frac{1}{2}} {\cal B} \mathcal{P}_{D}^{-\frac{1}{2}} \right) \\
&=& \left(\begin{array}{ccc}
{I} & {M^{T}} & {0} \\
{0} & {I} & {-N^{T}} \\
{0} & {N} & {0}
\end{array}\right)^{-1}
\left(\begin{array}{ccc}
{I} & {M^{T}} & {0} \\
{-M} & {0} & {-N^{T}} \\
{0} & {N} & {0}
\end{array}\right),
\end{eqnarray}
where $M=S^{-\frac{1}{2}} B A^{-\frac{1}{2}}$ and  $N=(CS^{-1}C^T)^{-\frac{1}{2}} C S^{-\frac{1}{2}}$. It straightforward to verify that $MM^T=I$, $NN^T=I$ and
\[
\hat{{\cal H}}=I+
\left(\begin{array}{ccc}
M^T(I-N^TN)M & {M^{T}(I-N^TN)} & {0} \\
{(N^TN-I)M} & {N^TN-I} & 0 \\
{NM} & {N} & {0}
\end{array}\right).
\]
Direct computation reveals that $(\hat{{\cal H}}-I)^2=0$. This shows that the minimal polynomial of $\hat{{\cal H}}$, as well as ${\cal H}$ is 2.
\end{proof}
\end{theorem}

\begin{remark}\label{remG}
	Theorem \ref{minpol} shows that the complete version of the GMRES method for solving the system  ${\cal P}^{-1} {\cal B}{\bf x}={\cal P}^{-1}\bf{b}$ will converge in two iterations in exact arithmetic.
\end{remark}
	
\section{Implementation of the preconditioner}\label{sec4}

In the implementation of the preconditioner ${\cal P}$ in a Krylov subspace method like GMRES,  in each iteration, a vector of the form
$v=\mathcal{P}^{-1}w$ should be computed.
To this end, all we need is to solve $\mathcal{P} v=w$ for $v$. If we set $v=(v_{1};v_{2};v_{3})$ and $w=(w_{1};w_{2};w_{3})$ in
which $v_{1}, w_{1} \in \mathbb{R}^{n},$ $v_{2}, w_{2} \in \mathbb{R}^{m}$ and $v_{3}, w_{3} \in \mathbb{R}^{l}$, then
we need to solve the system
\begin{equation*}
\left(\begin{array}{ccc}
{A} & {B^{T}} & {0} \\
{0} & {S} & {-C^{T}} \\
{0} & {C} & {0}
\end{array}\right)\left(\begin{array}{l}
{v_{1}} \\
{v_{2}} \\
{v_{3}}
\end{array}\right)=\left(\begin{array}{l}
{w_{1}} \\
{w_{2}} \\
{w_{3}}
\end{array}\right).
\end{equation*}
The following algorithm is given for solving the above linear system of equations.

\noindent\hrulefill

\noindent Algorithm 1: Computation  of $(v_1;v_2;v_3)=\mathcal{P} ^{-1}(w_1;w_2;w_3)$.

\noindent\hrulefill\\
\noindent 1: Set $t_{1}=w_{3}-C S^{-1} w_{2};$\\
2: Solve $(C S^{-1} C^{T})v_{3}=t_{1}$ using the Cholesky factorization of $C S^{-1} C^{T}$;\\
3: Set  $t_{2}=w_{2}+C^{T} v_{3};$\\
4: Solve $S v_{2}= t_{2}$ ;\\
5: Set $t_{3}=w_{1}-B^{T} v_{2};$\\
6: Solve $A v_{1}=t_{3}$ by the Cholesky factorization of $A.$

\noindent\hrulefill

\medskip

We end this section by pointing out to the choice of SPD matrix $S$. As seen, Remark \ref{remG} shows that $S=BA^{-1}B^T$ leads to an ideal case. However, by this choice, the resulting algorithm can be costly in general cases.
Basically, a preconditioner is called ``optimal", if the number of preconditioned iterations is independent of the size of the
problem and the amount of work per iteration scales
linearly with the size of the problem.  Notice that for our test problems,
total work (and, approximately, the corresponding CPU-time) should grow by
a factor of 4 each time the value of $p$ doubles.

In view of Remark \ref{remG} and the above discussions, in the numerical experiments,  we are particularly inspired to set $S=I$, with $I$ being the identity matrix or $S=\diag(B~  \diag(A)^{-1}  B^{T})$. For these choices, the proposed preconditioners, while not quite optimal, scales well with increasing the size of problem for our test examples. We also use $S=I$ or $S=\diag(B~  \diag(A)^{-1}  B^{T})$ while working with the preconditioners $\mathcal{P}_D$ and $\mathcal{P}_1$. In this work, we examine the exact versions of preconditioners in conjunction with complete version of GMRES.

In general cases, for approximating  $B A^{-1} B^{T}$ by $S$, similar to  \cite{Beik1,Beik2,Benzi2}, one can possibly avoid forming $B A^{-1} B^{T}$ and $CS^{-1}C^T$. Instead,  using a prescribed tolerance, few steps of the (P)CG method can be used for the actions of $A^{-1}$,  $(B A^{-1} B^{T})^{-1}$ and $(CS^{-1}C^T)^{-1}$. For this inexact implementation, the preconditioner should be used in conjunction with flexible GMRES (FGMRES). For some problems, we may have access to the sparse matrix $M$, being spectrally equivalent to $B A^{-1} B^{T}$.  In this case we can set $S=M$ and implement the preconditioner either exactly in conjunction with GMRES or inexactly in conjunction with FGMRES.

     \section{Numerical experiments}\label{sec5}
\begin{figure}[!tp]
      	\begin{center}
      		\includegraphics[height=3.5cm,width=3.5cm]{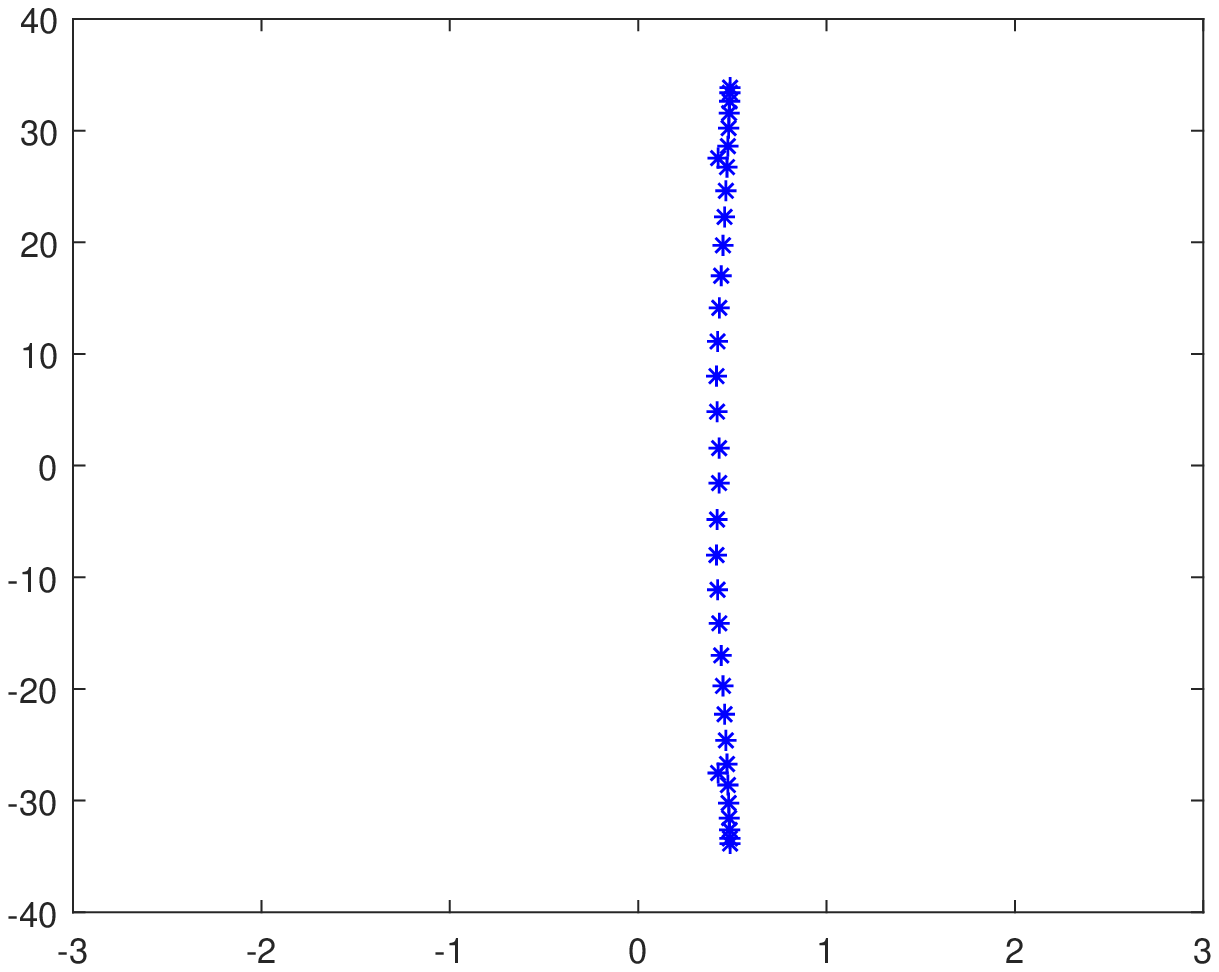}
      		\hspace*{0.25cm}
      		\includegraphics[height=3.5cm,width=3.5cm]{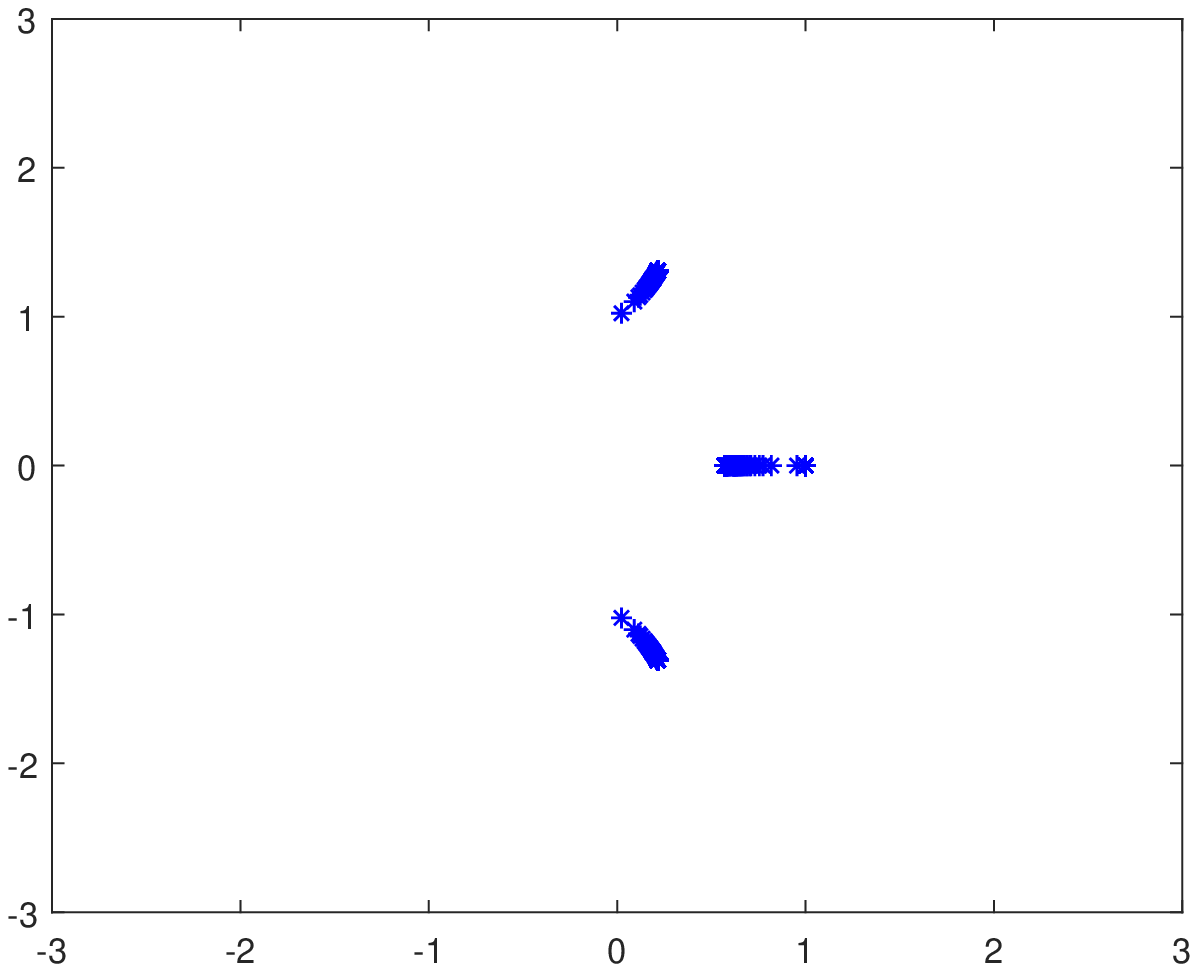}
      		\hspace*{0.25cm}		
      		\includegraphics[height=3.5cm,width=3.5cm]{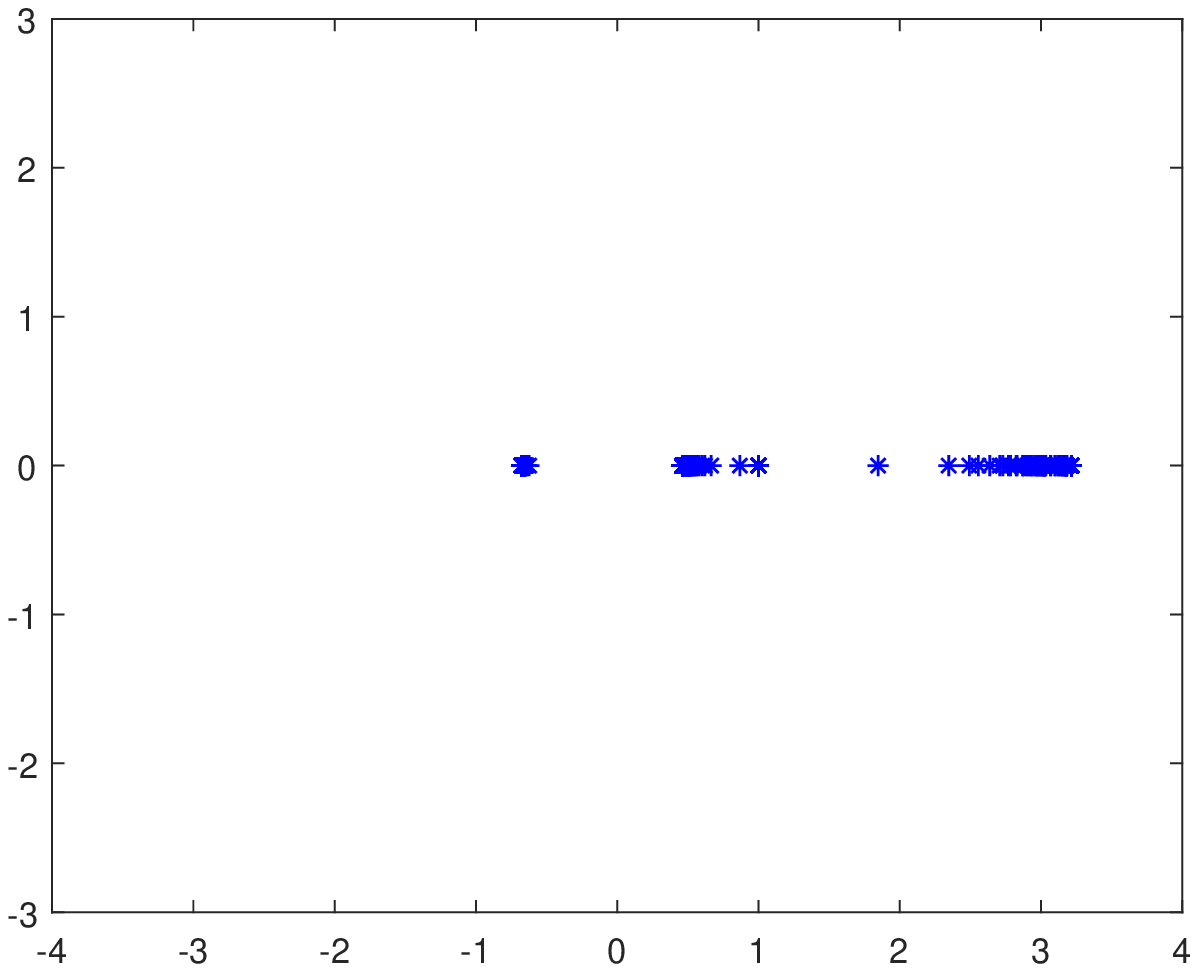}
      		\hspace*{0.25cm}
      		\includegraphics[height=3.5cm,width=3.5cm]{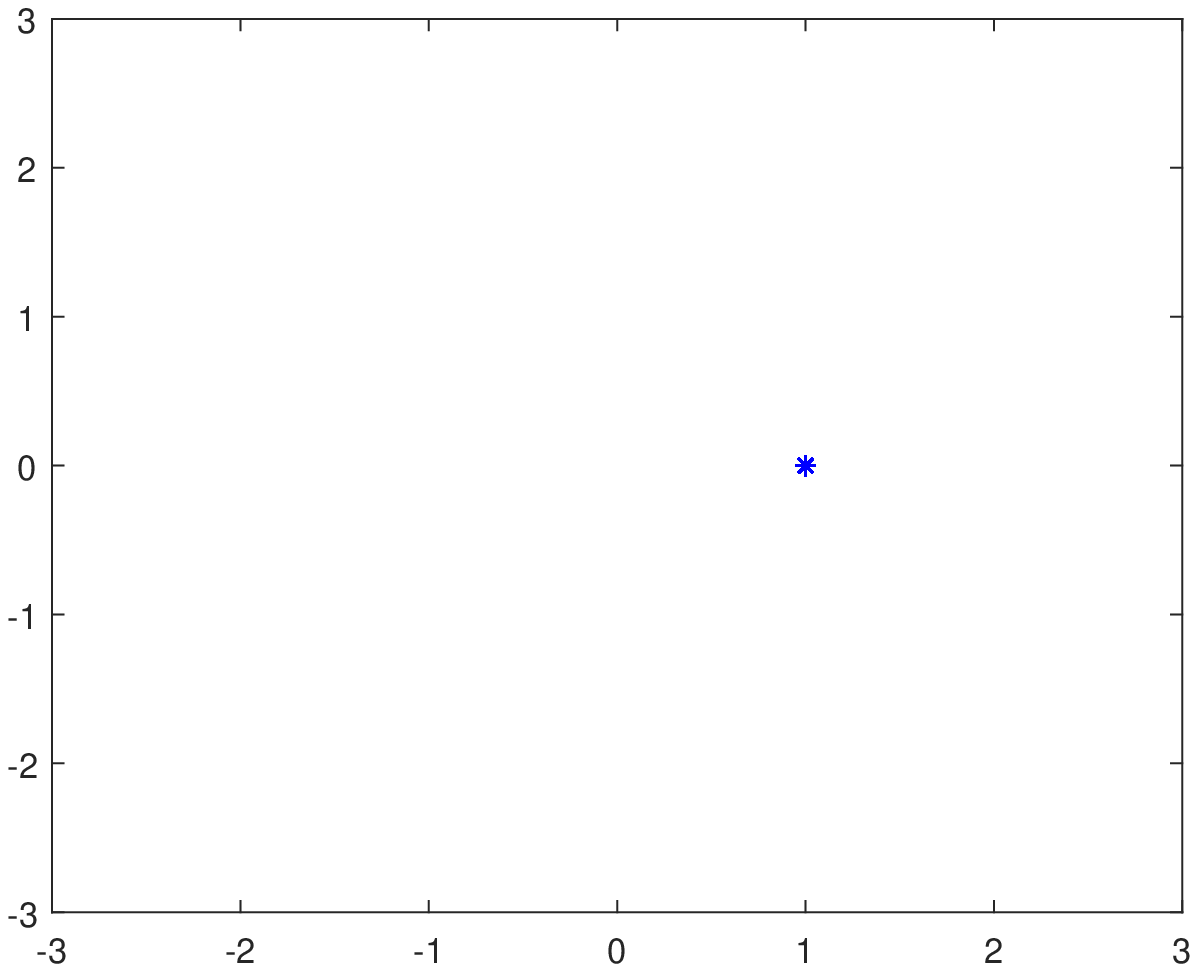}
      	\end{center}
      	\caption{{\small Eigenvalue distributions of $\mathcal{B},$ $\mathcal{P}_{D}^{-1} \mathcal{B},$ $\mathcal{P}_{1}^{-1} \mathcal{B}$ and $\mathcal{P}^{-1} \mathcal{B}$ (from the left to right) with $S=I$ and $p=16$ for Example \ref{ex1}}. \label{fg1} }	
      \end{figure}

     In this section, we numerically solve the  three-by-three saddle point problem \eqref{eq111} to examine the performance of proposed preconditioner in Section \ref{sec3}. In order to compare the performance of our preconditioner with the recently proposed ones in the literature, test problems are taken from \cite{A1,Huang-NumerAlgor,A2}. In all the test examples we use the complete version of GMRES method with right preconditioning. All runs were started from the initial zero vector and terminated once the current iteration (${\bf x}^{(k)}$) satisfies
     $$\frac{ \| {\bf b}-\mathcal{B} {\bf x} ^{(k)} \|}{\| {\bf b}\|}< 10^{-7},$$
     or the maximum number of iterations exceeds 5000.
     In all tests, the right-hand side vector ${\bf b}$ is set ${\bf b}=\mathcal{B}e$, where $e\in\Bbb{R}^{n+m+l}$ is vector of all ones.  Numerical results are presented in the tables in which ``IT" and ``CPU" denote the number of iterations and elapsed CPU times in second, respectively. A ``-" means that the method has not converged in the maximum number of iterations. To show the accuracy of the methods we also report the values
     \[
     Err=\frac{\|x^{(k)}- x^{*} \|} {\| x^{*} \|},
     \]
     in the tables, where $x^{*}$ stands for the exact solution of the system \eqref{eq111}. All runs were performed in  \textsc{Matlab} R2017a with a personal computer with 2.40 GHz central processing unit (Intel(R) Core(TM) i7-5500), 8 GB memory and Windows 10 operating system.
     \begin{example}{\rm \cite{A1,A2}} \label{ex1} \rm
      	Consider the saddle point problem \eqref{eq111} with
      	\begin{align*}
      	A=\left(\begin{array}{cc}
      	{I \otimes T+T \otimes I} & {0} \\
      	{0} & {I \otimes T+T \otimes I}
      	\end{array}\right) \in \mathbb{R}^{2 p^{2} \times 2 p^{2}},
      	\end{align*}
      	$
      	B=(I \otimes F \quad F \otimes I) \in \mathbb{R}^{p^{2} \times 2 p^{2}}$ and  $C=E \otimes F \in \mathbb{R}^{p^{2} \times p^{2}}
      	$
      	where
      	\begin{align*}
      	T=\frac{1}{h^{2}} \cdot \operatorname{tridiag}(-1,2,-1) \in \mathbb{R}^{p \times p}, \quad F=\frac{1}{h} \cdot \operatorname{tridiag}(0,1,-1) \in \mathbb{R}^{p \times p},
      	\end{align*}
      	and $
      	E=\operatorname{diag}\left(1, p+1, 2p+1, \ldots, p^{2}-p+1\right)
      $ in which $\otimes$ denotes Kronecker product and $h={1}/{(p+1)}$ the discretization meshsize.
      	
      	In this example, we set $S=I$, where $I$ is the identity matrix. Table \ref{tab11} shows the iteration counts and the elapsed CPU time for the GMRES method with the preconditioner  $\mathcal{P}_{D},$ ${\mathcal{P}_1}$ and $\mathcal{P}$.  To see the effectiveness of preconditioners, we have also reported the numerical results of the GMRES method without preconditioning.
      	Numerical results illustrate that the preconditioners can significantly reduce the number of iterations and elapsed CPU time of the GMRES method without preconditioning. As seen,  ${\cal P}$ is superior to the other examined preconditioners. An interesting observation which can be posed here is that the GMRES method with the preconditioner ${\cal P}$ gives the best accuracy among the preconditioners.  We also see that, $S=I$ presents a good approximation of the matrix $S=B A^{-1}B^{T}$.
      	
      	Fig. \ref{fg1} plots the eigenvalues of the matrices ${\cal B}$ ,  $\mathcal P_{D}^{-1} \mathcal{B},$  ${\mathcal P}_1^{-1} \mathcal{B}$  and $\mathcal P^{-1} \mathcal{B}$ for $p=16$ with $S=I$.  It is seen that the eigenvalues of $\mathcal P^{-1} \mathcal{B}$ are more clustered than the others.

      	\begin{table}[!t]
  		\centering
%      		\caption{Numerical results for Example \ref{ex1} with $S=\diag(B A^{-1} B^{T})$.\label{tab21}}\vspace{0.25cm}
%      		\begin{tabular}{|p{1.5cm}| p{1.5cm}|p{1.7cm}p{1.7cm} p{1.7cm}p{1.7cm}| }
%      			%\multicolumn{6}{l}{}\\
%      			\hline
%      			
%      			\multirow{3}{*}{Precon.}   & $p$          & 32     & 48   &  64   &  128 \\
%      			                           & $\mathbf{n}$ & 4096   & 9216 & 16384 & 65536  \\
%      			\hline\hline
%      			
%      			\multirow{3}{*}{$I$}& IT   & 3100     & -     & -    & -\\
%      			                    & CPU  & 89.40    & -     & -    & -\\
%      			                    & Err  & 8.86e-08 & -     & -    & -\\
%      			\hline
%      			\multirow{3}{*}{$\mathcal{P}_{D}$}&
%      			IT&38&39&40&44\\
%      			& CPU&0.41&1.66&5.76&111.58\\
%      			&Err&4.00e-06&5.07e-06&1.68e-05&3.33e-05\\
%      			\hline
%      			\multirow{3}{*}{$\mathcal{P}_1$}&
%      			IT&29&31&32&32\\
%      			& CPU&0.31&1.54&5.59&112.50\\
%      			&Err&9.11e-07&3.50e-07&6.07e-07&1.73e-05\\
%      			\hline
%      			\multirow{3}{*}{$\mathcal{P}$}&
%      			IT&2&2&2&2\\
%      			& CPU&0.26&1.42&5.37&110.78\\
%      			& Err&1.40e-12&4.72e-12&1.18e-11&8.75e-11\\
%      			\hline
%      		\end{tabular}\\[0.5cm]
      		
%      	\end{table}
%      		\begin{table}[]
      		\caption{Numerical results for Example \ref{ex1} with $S=I$.\label{tab11}}\vspace{0.25cm}
      		\begin{tabular}{|p{1.5cm}| p{1.5cm}|p{1.7cm}p{1.7cm} p{1.7cm}p{1.7cm}| }
      			%\multicolumn{6}{l}{}\\
      			\hline

      			\multirow{3}{*}{Precon.}   & $p$          &   64   &  128  &  256   &  512\\
      		                             & $\mathbf{n}$   &  16384 & 65536 & 262144 &  1048576 \\
      			\hline\hline

      			\multirow{3}{*}{$I$}&
      			IT&-&-&-&-\\
      			& CPU&-&-&-&-\\
      			&Err&-&-&-&-\\
      			\hline
      			\multirow{3}{*}{$\mathcal{P}_{D}$}&
      			IT&36&39&41&47\\
      			& CPU&0.38&1.90&13.73&114.90\\
      			&Err&1.46e-05&1.33e-05&1.08e-04&1.13e-04\\
      			\hline
      			\multirow{3}{*}{$\mathcal{P}_1$}&
      			IT&28&30&30&32\\
      			& CPU&0.26&1.24&10.99&86.42\\
      			&Err&2.08e-06&6.50e-06&2.95e-05&5.79e-05\\
      			\hline
      			\multirow{3}{*}{$\mathcal{P}$}&
      			IT&2&2&2&6\\
      			& CPU&0.06&0.39&2.61&34.05\\
      			& Err&1.16e-11&6.50e-11&6.84e-10&5.02e-09\\
      			\hline
      		\end{tabular}
      		\end{table}
       \end{example}

 %================================ Example 2 %%%%%%%%%%%%%%%%%%%%%%%%%%%%%%%%%%%%%%

     \begin{example} {\rm \cite{A1,A2}}\label{xe1} \rm
    	Consider the three-by-three block saddle point problem \eqref{eq1} for which
    	\begin{align*}
    	A=\operatorname{bldiag}\left(2 W^{T} W+D_{1}, D_{2}, D_{3}\right) \in \mathbb{R}^{n \times n},
    	\end{align*}
    	is a block-diagonal matrix,
    	\begin{align*}
    	B=\left[E,-I_{2 \widetilde{p}}, I_{2 \widetilde{p}}\right] \in \mathbb{R}^{m \times n} \text { and } \quad C=E^{T} \in \mathbb{R}^{l \times m},
    	\end{align*}
    	are both full row-rank matrices where $\tilde{p}=p^{2},$
    	$\hat{p}=p(p+1);$  $D_{1}=I_{\widehat{p}}$ is an identity matrix; $D_{i}=\operatorname{diag}(d_{j}^{(i)}) \in \mathbb{R}^{2\tilde{p} \times 2\tilde{p}},$ $i=2,3,$ are diagonal matrices, with
    	\begin{eqnarray*}
    	d_{j}^{(2)}&\h=\h&\left\{\begin{array}{ll}
    	{1,} & {\text { for } \quad 1 \leq j \leq \tilde{p}}, \\
    	{10^{-5}(j-\tilde{p})^{2}},  & \text { for }   {\tilde{p}+1 \leq j \leq 2 \tilde{p}},
    	\end{array}\right. \\
    	%\end{align*}
    	%\begin{align*}
    	d_{j}^{(3)} &\h=\h&10^{-5}(j+\tilde{p})^{2}
    	\text { for } 1 \leq j \leq 2 \tilde{p},
    	\end{eqnarray*}
    	and
    	\begin{align*}
    	E=\left(\begin{array}{cc}
    	{\widehat{E} \otimes I_{p}} \\
    	{I_{p} \otimes \widehat{E}}
    	\end{array}\right), \quad \widehat{E}=\left(\begin{array}{ccccc}
    	{2} & {-1} & {} & {} & {}\\
    	{} & {2} & {-1} & {} & {} \\
    	{} & {} & {\ddots} & {\ddots} & {}\\
    	{} & {} & {} & {2} & {-1}
    	\end{array}\right) \in \mathbb{R}^{p \times(p+1)}.
    	\end{align*}
    	Moreover, $W=v v^{T} \in \mathbb{R}^{\hat{p} \times \hat{p}}$, where  $v\in\mathbb{R}^{\hat{p}}$ is an arbitrary vector.
    		According to the above definitions, we have $n=\hat{p}+4\tilde{p}$, $m=2\tilde{p}$ and $l=\hat{p}$.
    		
    	We consider two choices for the vector $v$. In the first choice, the $i$th entry of the vector $v$ is set to be $v_{i}=e^{-2(i / 3)^{2}}$, $i=1,2,\ldots,l$, and in the second one the vector $v$ is set to be a random sparse vector of order $l$ with
    	approximately $0.05l$ uniformly distributed nonzero entries (such a vector can be generated using the ``\verb'sprand'" command of \textsc{Matlab}). 	
    		
     	For both of the choices we set $S=I$. Numerical results for the first choice are presented  in Table \ref{tab3} and for the second choice   in Table \ref{tab4}.  All the other notations are as the previous example.  As seen, the proposed preconditioner outperforms the others in terms of the iteration counts, the elapsed CPU time and the accuracy of computed solution. Fig. \ref{fig2} and  Fig. \ref{fig3} display the eigenvalue distribution of the original coefficient matrix, $\mathcal{P}_{D}^{-1} \mathcal{B},$ ${\cal P}_1^{-1} \mathcal{B}$ and $\mathcal{P}^{-1} \mathcal{B}$ for $S=I$ and $p=16$ for the two  choices, respectively. As observed,  eigenvalues of     	$\mathcal{P}^{-1} \mathcal{B}$ are more clustered around the point $(1,0)$ than the others.

    	\begin{table}[!t]
    		\centering
    		\caption{Numerical results for Example \ref{xe1} for the first choice with $S=I$.\label{tab3}}\vspace{0.25cm}
      		\begin{tabular}{|p{1.5cm}| p{1.5cm}|p{1.7cm}p{1.7cm} p{1.7cm}p{1.7cm}| }
      			%\multicolumn{6}{l}{}\\
      			\hline
      			
      			\multirow{3}{*}{Precon.}   & $p$          & 32     & 48   &  64   &  128 \\
      			                           & $\mathbf{n}$ & 8256   & 9216 & 32896 &   131328  \\
      			\hline\hline

    			\multirow{3}{*}{$I$}&
    			IT&557&1180&1815&2128\\
    			& CPU&5.22&38.50&133.15&209.50\\
    			&Err&5.22-06&5.67e-05&1.84e-04&6.61e-04\\
    			\hline
    			\multirow{3}{*}{$\mathcal{P}_{D}$}&
    			IT&348&314&284&197\\
    			& CPU&10.97&13.51&14.84&21.75\\
    			&Err&4.87e-06&2.34e-05&7.51e-05&1.14e-03\\
    			\hline
    			\multirow{3}{*}{${\mathcal{P}}_1$}&
    			IT&171&159&144&103\\
    			& CPU&2.73&3.74&4.09&7.24\\
    			&Err&3.48e-06&2.19e-5&7.43e-05&1.07e-03\\
    			\hline
    			\multirow{3}{*}{$\mathcal{P}$}&
    			IT&2&2&2&2\\
    			& CPU&0.08&0.17&0.06&0.40\\
    			& Err&5.64e-09&1.00e-08&2.06e-08&1.82e-08\\
    			\hline
    	\end{tabular}\vspace{0.5cm}
%\end{table}
%\begin{table}[]
    	\caption{Numerical results for Example \ref{xe1} for the second choice with $S=I$.\label{tab4}}\vspace{0.25cm}
    	     \begin{tabular}{|p{1.5cm}| p{1.5cm}|p{1.7cm}p{1.7cm} p{1.7cm}p{1.7cm}| }
    	      %\multicolumn{6}{l}{}\\
    	      \hline
    	      \multirow{3}{*}{Precon.}   & $p$            &   64   &  128   &  256   &  512\\
    	                                 & $\mathbf{n}$   &  32896 & 131328 & 524800 & 2098176  \\
    	      			\hline\hline
    	
    			\multirow{3}{*}{$I$} & IT  & -     &-&-&-\\
    			                     & CPU & -  &-&-&-\\
    			                     & Err & - &-&-&-\\
    			\hline
    			\multirow{3}{*}{$\mathcal{P}_{D}$}&
    			IT&279&193&125&119\\
    			& CPU&14.36&21.92&93.29&201.56\\
    			&Err&1.31e-04&2.19e-03&2.19e-02&3.72e-02\\
    			\hline
    			\multirow{3}{*}{$\mathcal{P}_1$}&
    			IT&143&103&70&59\\
    			& CPU&4.08&8.02&56.08&149.11\\
    			&Err&1.30e-04&1.94e-03&2.17e-04&2.12e-03\\
    			\hline
    			\multirow{3}{*}{$\mathcal{P}$}&
    			IT&2&2&2&4\\
    			& CPU&0.08&0.78&30.92&66.42\\
    			& Err&1.33e-09&4.09e-09&3.20e-09&2.41e-09\\
    			\hline
    		\end{tabular}
    		\label{tab4}
    	\end{table}

    	\begin{figure}[!t]
     		\begin{center}
    			\includegraphics[height=3.5cm,width=3.5cm]{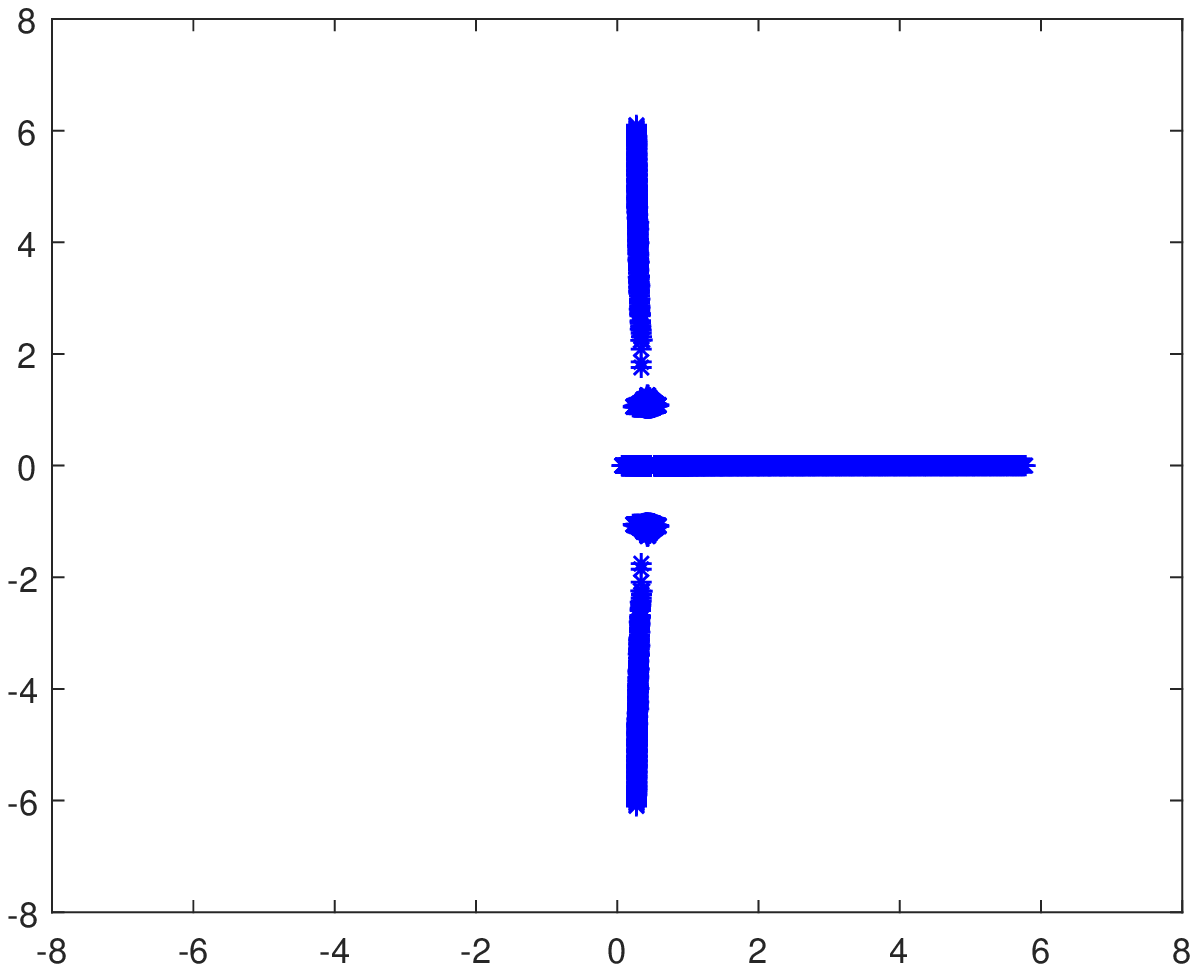}
    			\hspace*{0.25cm}
    			\includegraphics[height=3.5cm,width=3.5cm]{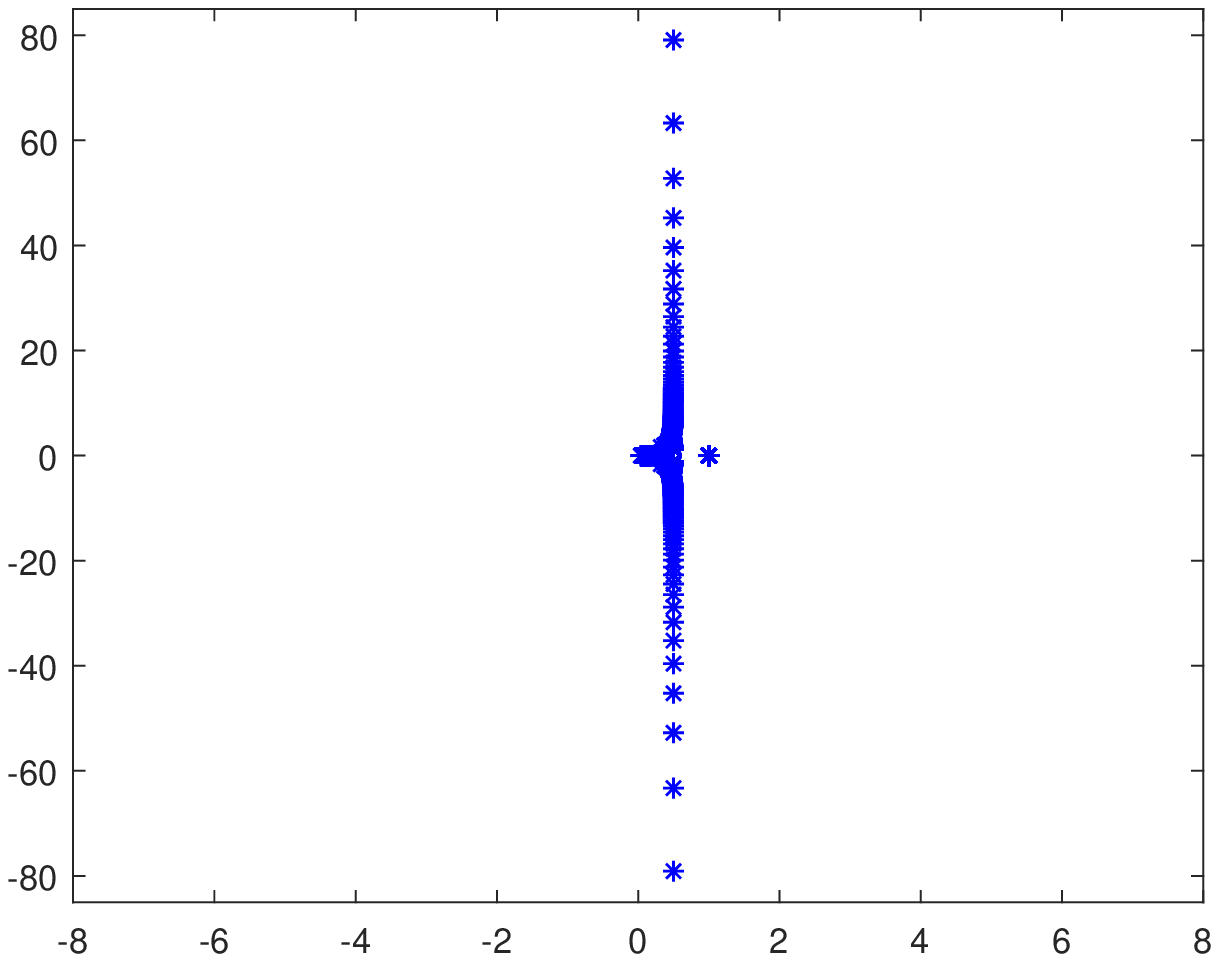}
    			\hspace*{0.25cm}
    			\includegraphics[height=3.5cm,width=3.5cm]{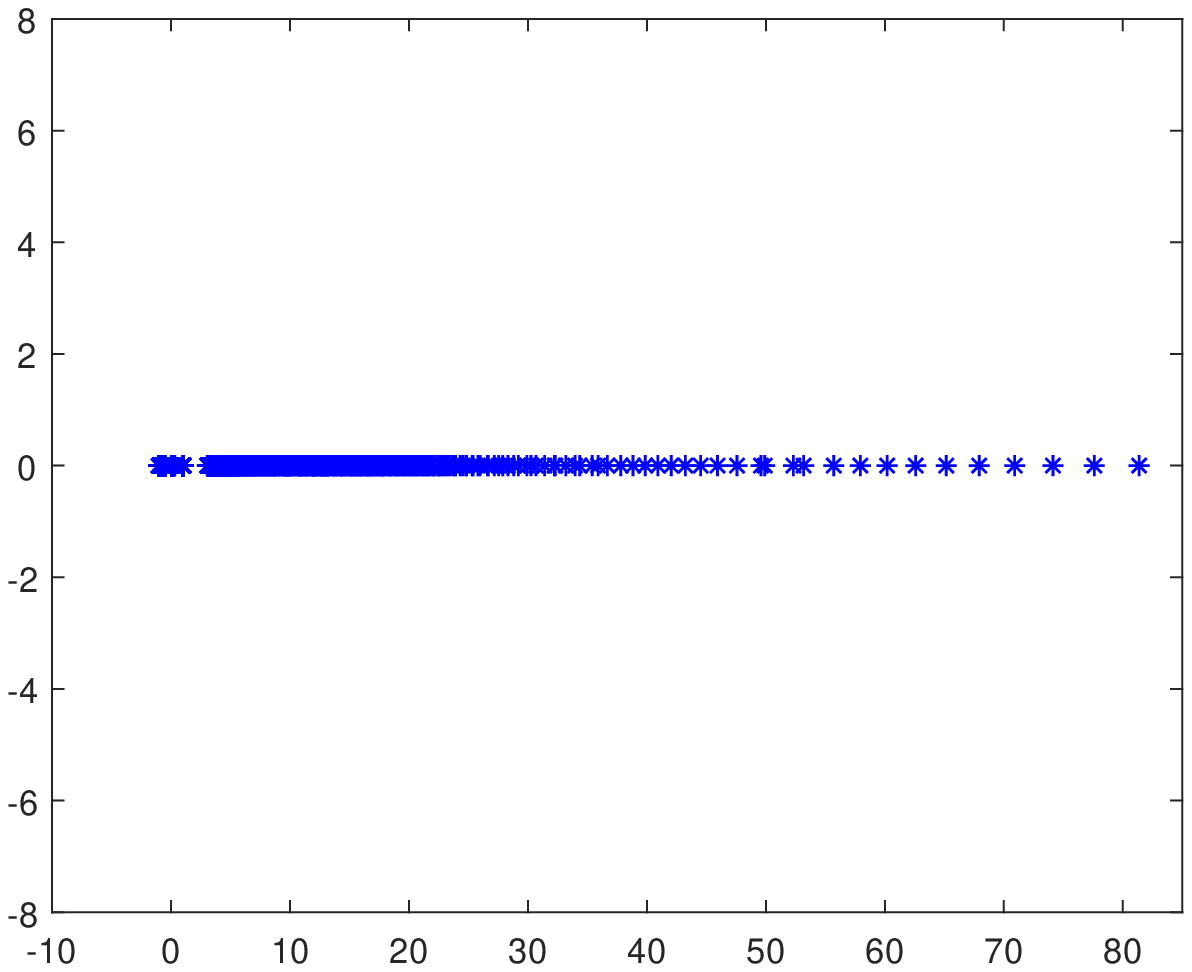}
    			\hspace*{0.25cm}
    			\includegraphics[height=3.5cm,width=3.5cm]{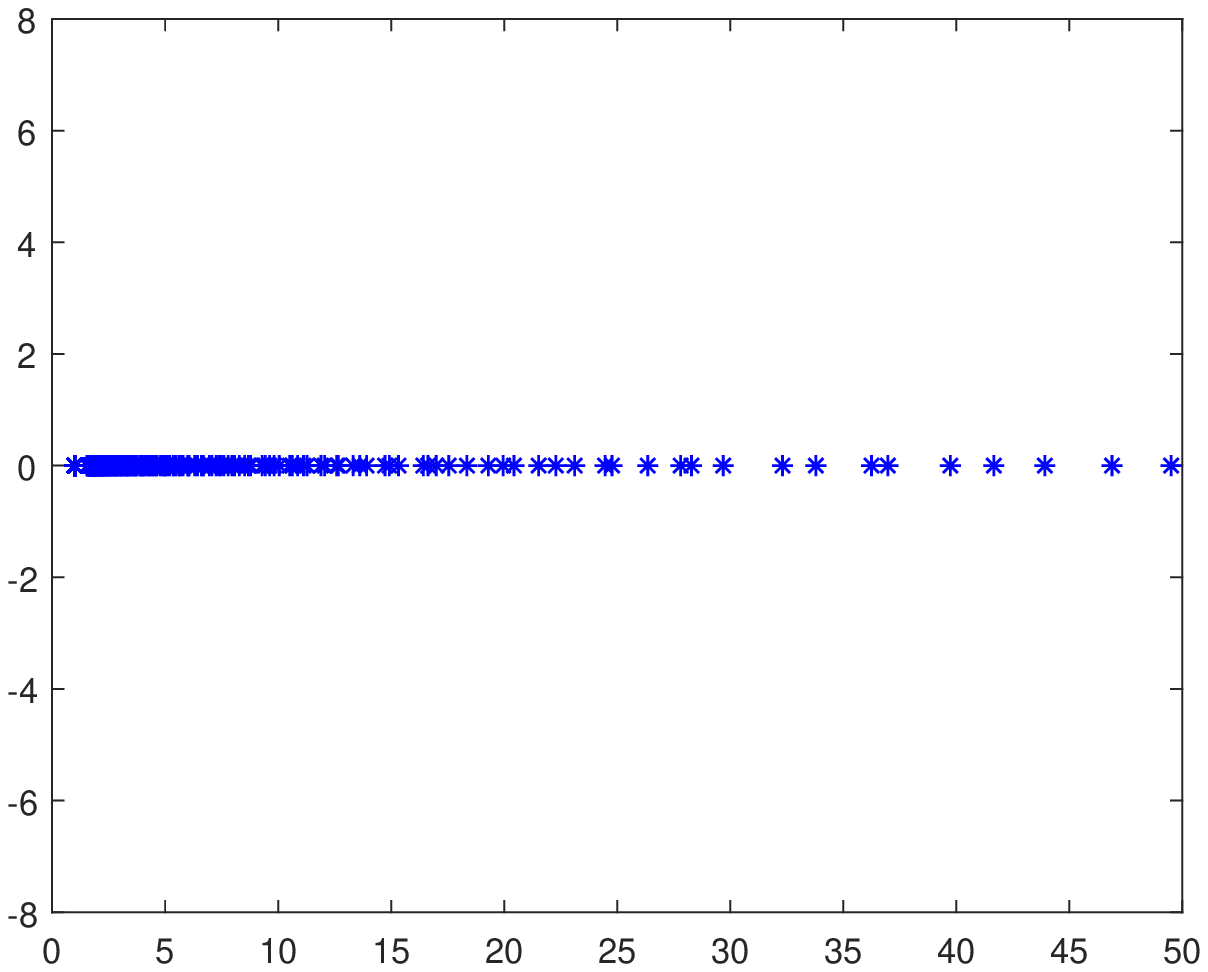}
    		\end{center}
    		\caption{{\small Eigenvalue distributions of $\mathcal{B},$ ${\mathcal{P}_{D}}^{-1} \mathcal{B}, $ ${\mathcal{P}}_1^{-1} \mathcal{B}$
    				 and ${\mathcal{P}}^{-1} \mathcal{B}$ for the first choice with \textcolor{blue}{$S=I$} and $p=16$ for Example \ref{xe1}}. \label{fig2}}
%    	\end{figure}
%
%    \begin{figure}[!tp]
   	\begin{center}
    		\includegraphics[height=3.5cm,width=3.5cm]{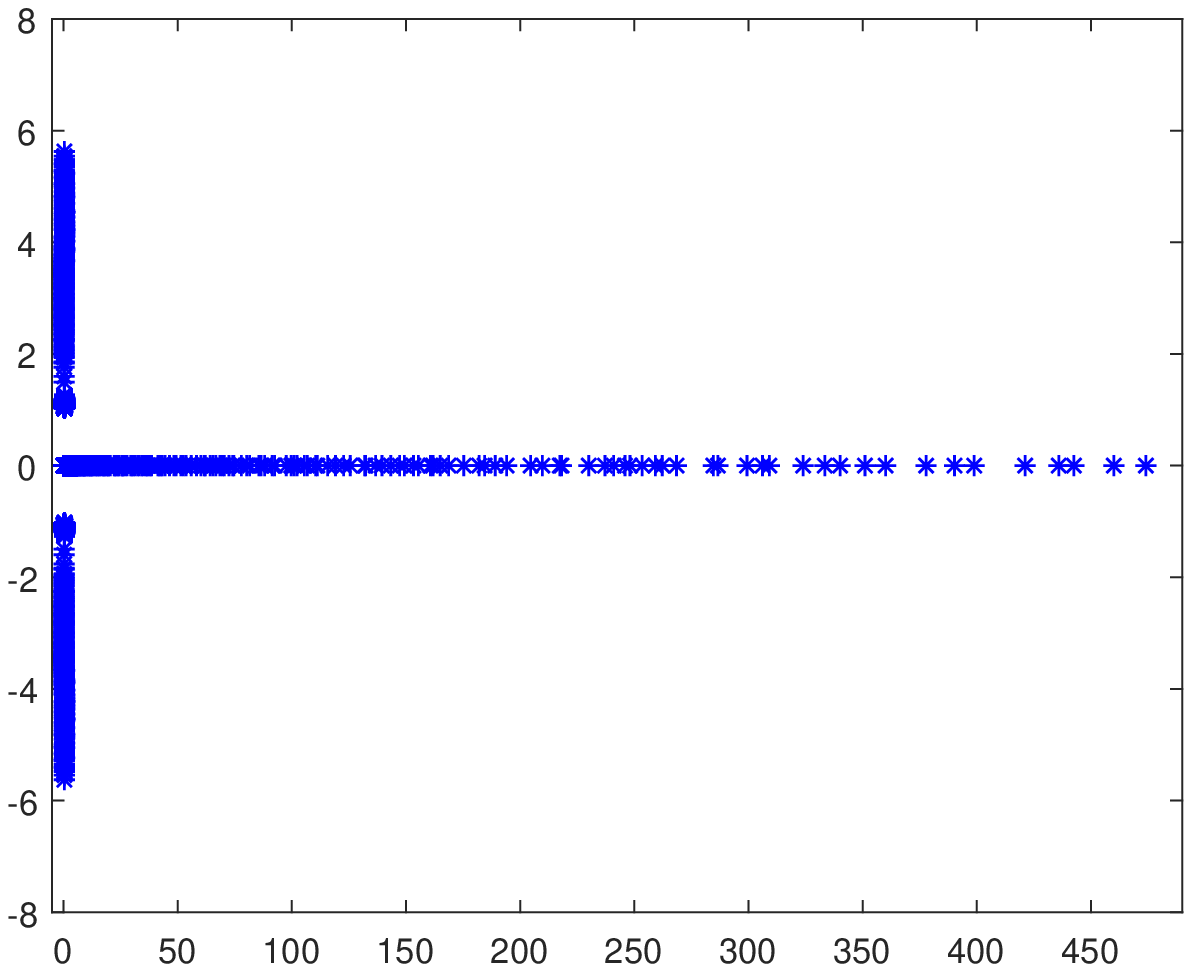}
    		\hspace*{0.25cm}
    		\includegraphics[height=3.5cm,width=3.5cm]{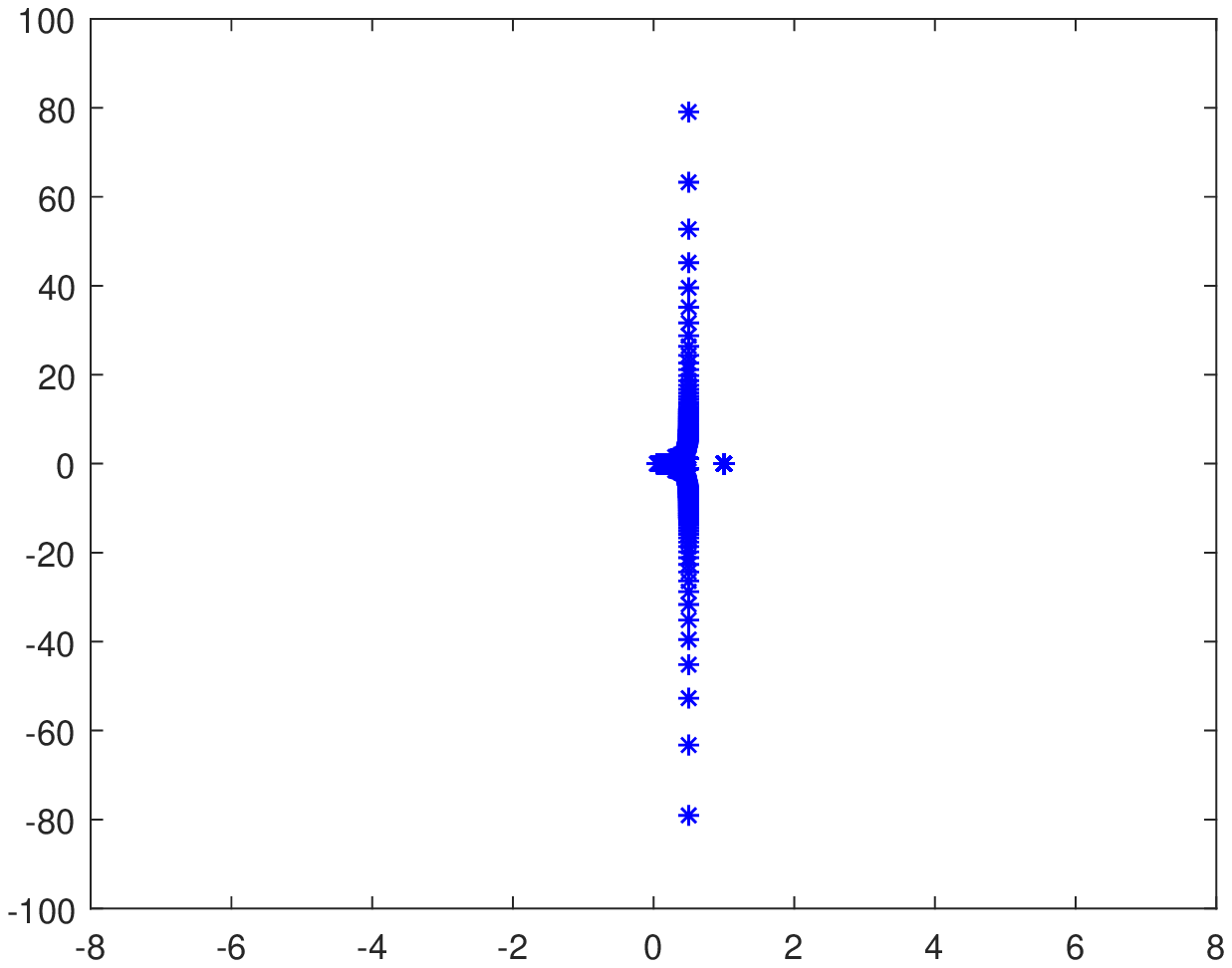}
    		\hspace*{0.25cm}
    		\includegraphics[height=3.5cm,width=3.5cm]{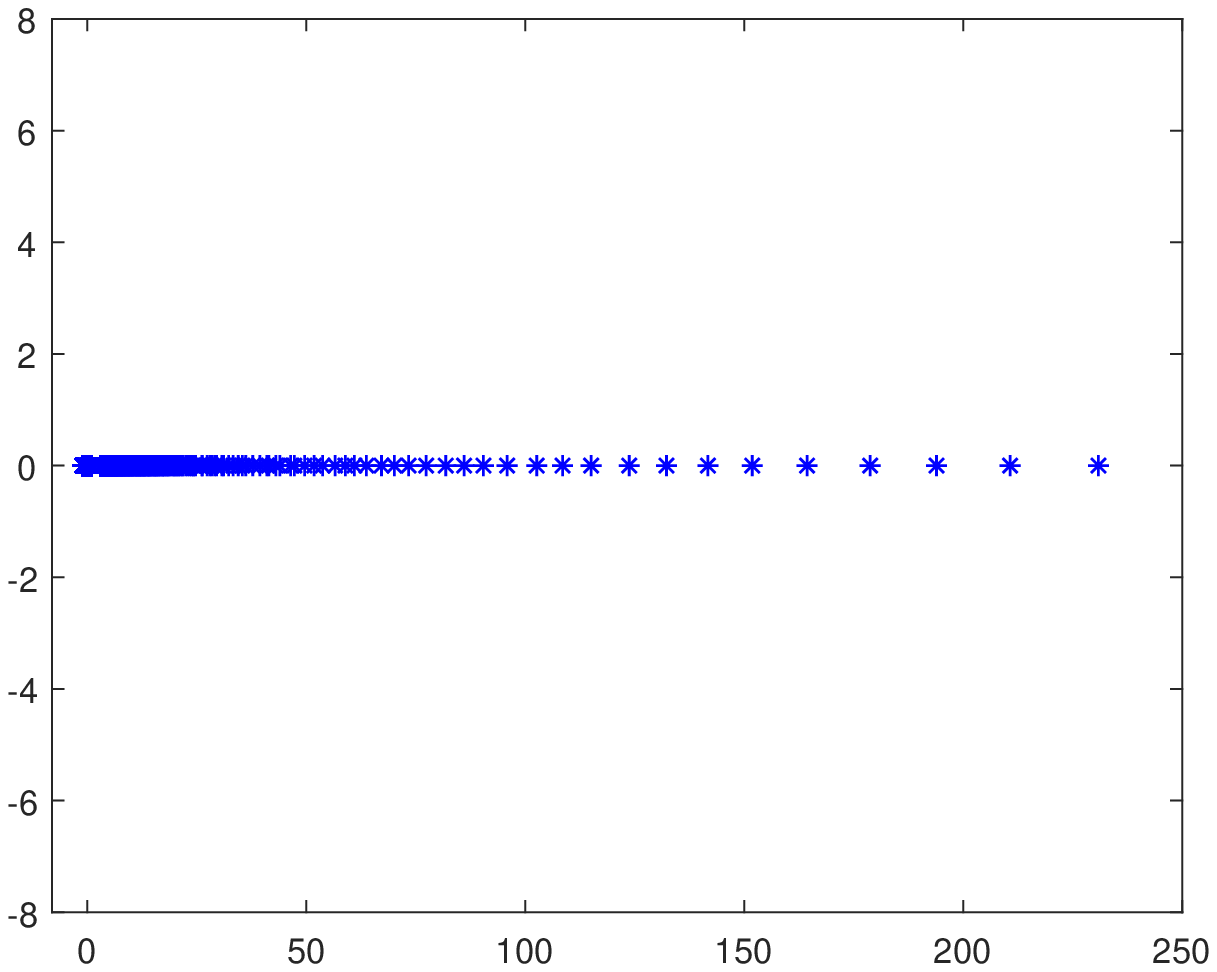}
    		\hspace*{0.25cm}
    		\includegraphics[height=3.5cm,width=3.5cm]{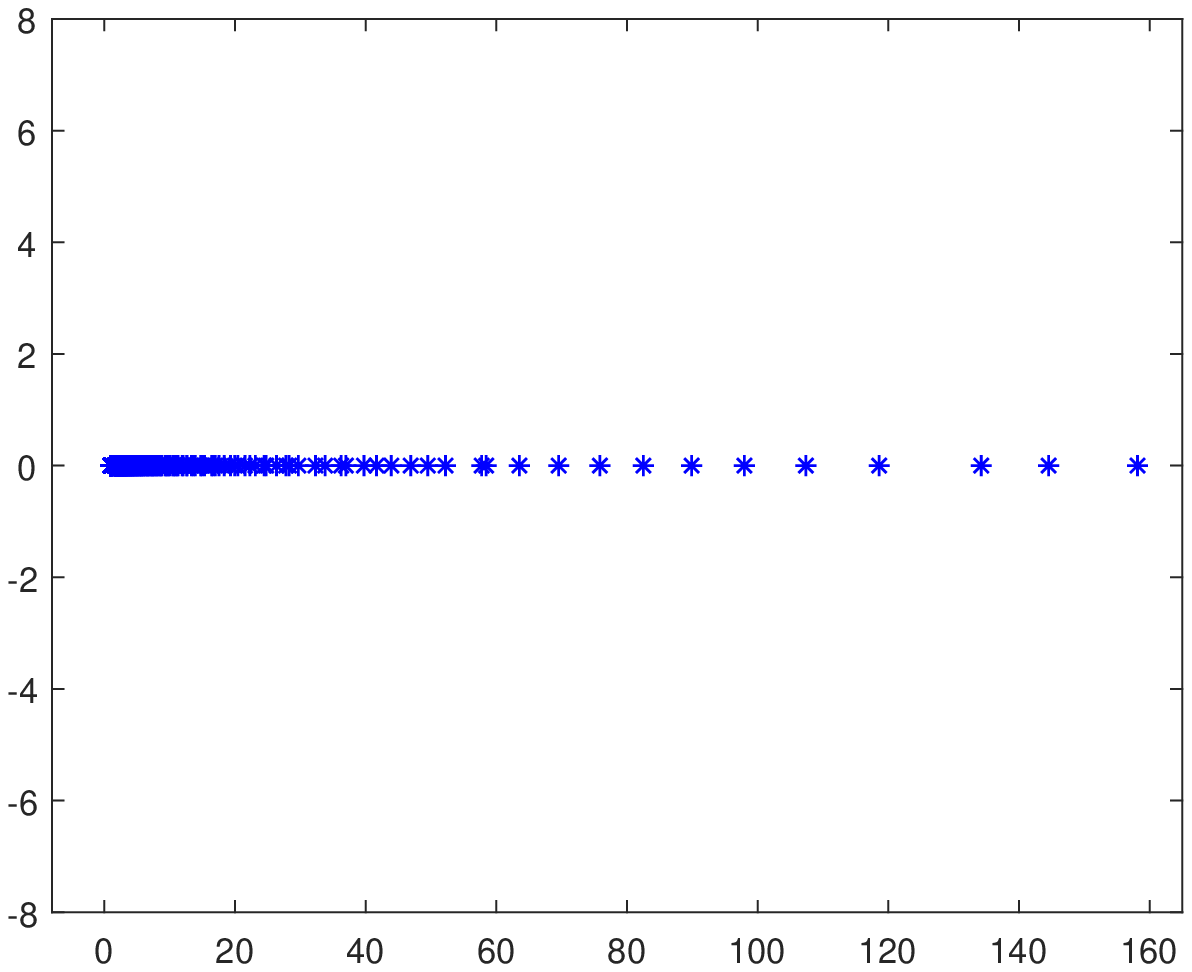}
    	\end{center}
    	\caption{{\small Eigenvalue distributions of $\mathcal{B},$ ${\mathcal{P}_{D}}^{-1} \mathcal{B}, $ ${\mathcal{P}}_1^{-1} \mathcal{B}$
    			and ${\mathcal{P}}^{-1} \mathcal{B}$ for the second choice with \textcolor{blue}{$S=I$} and $p=16$ for Example \ref{xe1}}. \label{fig3}}
    \end{figure}
    \end{example}
\begin{example}\label{ex3}\rm
		We consider the three-by-three block saddle point problem \eqref{eq111} with (see \cite{Huang-NumerAlgor,A2})
		\begin{eqnarray}
		&\h\h& \min_{x\in\Bbb{R}^n,y\in\Bbb{R}^l} \frac{1}{2}x^TAx+r^Tx+q^Ty \label{CUTER}\\
		\nonumber &\h\h&~~~s.t.:~~~ Bx+C^Ty=b,
		\end{eqnarray}
		where $r\in\Bbb{R}^n$ and $q\in\Bbb{R}^l$. To solve the above problem we define the Lagrange function
		\[
		L(x, y, \lambda) =\frac{1}{2}x^TAx+r^Tx+q^Ty+\lambda^T(Bx+C^Ty-b),
		\]
		where the vector $\lambda\in\Bbb{R}^{m}$ is the Lagrange multiplier. Then the Karush-Kuhn-Tucker necessary conditions of \eqref{CUTER} are as follows (see \cite{Bertsekas})
		\[
		\nabla_{x} L(x, y, \lambda)=0,\quad \nabla_{y} L(x, y, \lambda)=0 \quad \text{and} \quad \nabla_{\lambda} L(x, y, \lambda)=0.
		\]
		It is easy to see that these equations give a system of linear equations of the form \eqref{eq1}. In this example, we have chosen the matrices $A$, $B$ and $C$ from the CUTEr collection \cite{CUTErCol}. To do so, we have selected four matrices. In this example, we  set $S=\diag(B~  \diag(A)^{-1}  B^{T})$ (see \cite{Bergamaschi}). Numerical results are presented in Table \ref{tab555}. As we see the proposed preconditioner outperforms the others from the iteration counts, elapsed CPU time and accuracy of the computed solution point of view.

			\begin{table}[H]
			\centering	
			\caption{Numerical results for Example \ref{ex3} with $S=\diag(B~  \diag(A)^{-1}  B^{T})$. \label{tab555}}\vspace{0.25cm}

            \begin{tabular}{                  | p{1.5cm}| p{1.5cm}|p{2.2cm} p{2.2cm}p{2.2cm}p{2.2cm}|}\hline
                           \multirow{3}{*}{Precon.}& Matrix        & AUG2DC    & AUG3DC     & LISWET12   &   YAO      \\
                                                   & $\mathbf{n}$  & 50400     & 8746       & 30004      &  6004      \\ \hline\hline
                              \multirow{3}{*}{$I$} & IT            & 94        & 99         &  92        &  99        \\
                                                   & CPU           & 2.12      & 0.93       &  1.52      &  0.67     \\
				                                   & Err           & 4.82e-07  & 2.37e-07   &  5.80e-07  &  5.93e-07 \\ \hline
                \multirow{3}{*}{$\mathcal{P}_{D}$} & IT            & 101       & 136        &  52        &  57        \\
				                                   & CPU           & 2.62      & 1.85       &  0.66      &  0.31     \\
				                                   & Err           & 3.41e-07  & 1.84e-07   &  4.40e-07  &  3.11e-07 \\ \hline
				\multirow{3}{*}{${\mathcal{P}}_1$} & IT            & 55        & 80         & 34         &  37       \\
				                                   & CPU           & 0.92      & 0.65       & 0.32       &  0.14     \\
				                                   & Err           & 3.68e-07  & 1.90e-7    & 3.88e-07   &  3.10e-09 \\ \hline
				\multirow{3}{*}{$\mathcal{P}$}     & IT            & 22        & 29         &  4         &  4        \\
				                                   & CPU           & 0.30      & 0.14       &  0.07      &  0.04     \\
				                                   & Err           & 1.33e-07  & 1.14e-7    & 1.68e-14   &  1.69e-14 \\
				\hline
			\end{tabular}
		\end{table}
	
\end{example}

\section{Conclusions} \label{sec6}
A new stationary iterative method was constructed for solving a class of three-by-three block saddle point problems. We analyzed the convergence properties of the elaborated stationary method. We further examined the performance of induced preconditioner from the proposed method. More precisely, the eigenvalue distribution of the preconditioned matrix was studied. Our numerical tests illustrated that the proposed preconditioner is more effective than the other tested preconditioners in the literature.

\section*{Acknowledgments}
 The work of the second author is partially supported by University of Guilan.

%\end{thebibliography}


\begin{thebibliography}{99}

\bibitem{Assous} F. Assous, P. Degond, E. Heintze, P.A. Raviart,  J. Segre, On a finite-element method for solving the three-dimensional Maxwell equations, J. Comput. Phys. 109 (1993) 222--237.

\bibitem{Beik1} F. P. A. Beik and M. Benzi, Block preconditioners for saddle point systems arising from liquid crystal directors modeling, CALCOLO 55 (2018)  29.

\bibitem{Beik2} F. P. A. Beik and M. Benzi, Iterative methods for double saddle point systems, SIAM J. Matrix Anal. Appl. 39 (2018)  902--921.

\bibitem{Benzi2} M. Benzi and F. P. A. Beik, Uzawa-type and augmented lagrangian methods for double saddle point systems, Structur Matrices in Numerical Linear Algebra (Prof. Dario Andrea Bini, Prof. Fabio Di Benedetto, Prof. Eugene Tyrtyshnikov, and Prof. Marc Van Barel, eds.), Springer International Publishing, 2019.

\bibitem{Benzi1} M. Benzi, G. H. Golub and J. Liesen, Numerical Solution of Saddle Point Problems, Acta Numer. 14 (2005)  1--137.

\bibitem{BenziSIAM} M. Benzi, G.H. Golub, A preconditioner for generalized saddle point problems, SIAM J. Matrix Anal. Appl. 26 (2004) 20-41.

\bibitem{A11}	M. Benzi, M.K. Ng, Q. Niu, Z. Wang, A relaxed dimensional factorization preconditioner for the incompressible Navier-Stokes equations, J. Comput. Phys. 230 (2011) 6185--6202.

\bibitem{Bergamaschi} L. Bergamaschi, J. Gondzio, G. Zilli, Preconditioning indefinite systems in interior point
methods for optimization, Comput. Optim. Appl. 28 (2004) 149--171.

\bibitem{Bertsekas} D. P. Bertsekas,  Nonlinear Programming, 2nd Ed., Athena Scientific, 1999.	

\bibitem{CAOAML} Y. Cao, Shift-splitting preconditioners for a class of block three-by-three saddle point problems, Appl. Math. Lett. 96 (2019) 40--46.
 	
\bibitem{A12} Z.-H.	 Cao, Positive stable block triangular preconditioners for symmetric saddle point problems, Appl. Numer. Math. 57 (2007) 899--910.	

 \bibitem{A4} Z.-M. Chen, Q. Du,  J. Zou, Finite element methods with matching and nonmatching meshes for	Maxwell equations with discontinuous coefficients, SIAM J. Numer Anal. 37 (2000) 1542--1570.	
	
\bibitem{A14}  C.-R. Chen, C.-F. Ma, A generalized shift-splitting preconditioner for singular saddle point problems, Appl. Math. Comput. 269 (2015) 947--955.		 

\bibitem{A6} Z.-M. Chen, Q. Du, J. Zou, Finite element methods with matching and nonmatching meshes for Maxwell equations with discontinuous coefficients, SIAM J. Numer. Anal. 37 (2000) 1542--1570.

\bibitem{A5} P. Ciarlet, J. Zou, Finite element convergence for the Darwin model to Maxwell's equations,
RAIRO Math. Modelling Numer. Anal. 31 (1997) 213--249.


\bibitem{A13} H.C.	Elman, D.J. Silvester, A.J. Wathen,  Performance and analysis of saddle point preconditioners
	for the discrete steady-state Navier-Stokes equations, Numer. Math. 90 (2002) 665--688.
	
	
\bibitem{CUTErCol}  N. I. M. Gould, D. Orban, P. L. Toint, CUTEr and SifDec, a constrained and unconstrained testing
	environment, revisited, ACM Trans. Math. Softw. 29 (2003) 373–394.

\bibitem{A8}  D.R.	Han, X.M. Yuan, Local linear convergence of the alternating direction method of multipliers	for quadratic programs, SIAM J. Numer. Anal. 51 (2013) 3446--3457.
	
\bibitem{Horn1} R.A. Horn,  C.R. Johnson, Matrix Analysis, Cambridge University Press, Cambridge, UK, 1985.

\bibitem{A1}  N. Huang, C.-F. Ma, Spectral analysis of the preconditioned system for the 3 $\times$ 3 block saddle point problem, Numer. Algor. 81 (2019) 421--444.

\bibitem{Huang-NumerAlgor} N. Huang, Variable parameter Uzawa method for solving a class of block three-by-three saddle point problems, Numer. Algor., 2020, \url{https://doi.org/10.1007/s11075-019-00863-y}.

\bibitem{Huang-NLWA} N. Huang, Y.‐H. Dai,  Q. Hu, Uzawa methods for a class of block three-by-three saddle‐point problems,
Numer Linear Algebra Appl. 26 (2019) e2265.

\bibitem{A15}	Y.-F. Ke, C.-F. Ma, The parameterized preconditioner for the generalized saddle point problems from the incompressible Navier-Stokes equations, J. Comput. Appl. Math. 37 (2018) 3385--3398.

\bibitem{A9}  Y. Saad, Iterative Methods for Sparse Linear Systems, Second Edition, Society for Industrial and Applied Mathematics, Philadelphia, 2003.

\bibitem{Salkuyeh} D.K. Salkuyeh, M. Masoudi, A new relaxed HSS preconditioner for saddle point problems, Numer. Algor. 74 (2017) 781--795.


\bibitem{A2} X. Xie, H.-B. Li, A note on preconditioning for the $3 \times 3$ block saddle point problem, Comput. Math. Appl.
79 (2020) 3289--3296.

\bibitem{A7} J.-Y. Yuan, Numerical methods for generalized least squares problems, J. Comput. Appl. Math. 66 (1996) 571--584.

\bibitem{Zhang} F. Zhang,  Q. Zhang, Eigenvalue inequalities for matrix product, IEEE Trans. Automat. Control 51 (2006) 1506--1509.

\end{thebibliography}
\end{document}